\documentclass{aptpub}
\authornames{P.S.Griffin and R.A. Maller}
\shorttitle{Stability of the Exit Time}

\newcommand{\topru}{\stackrel{\mathrm{P^{(u)}}}{\longrightarrow}}
\newcommand{\toprstar}{\stackrel{\mathrm{P^*}}{\longrightarrow}}
\newcommand{\todrstar}{\stackrel{\mathrm{D^*}}{\longrightarrow}}

\newcommand{\be}{\begin{equation}}
\newcommand{\ee}{\end{equation}}
\newcommand{\ben}{\begin{equation*}}
\newcommand{\een}{\end{equation*}}
\newcommand{\ba}{\begin{aligned}}
\newcommand{\ea}{\end{aligned}}
\newcommand{\veps}{\varepsilon}
\newcommand{\rmd}{{\rm d}}
\newcommand{\rmi}{{\rm i}}

\newcommand{\dto}{\downarrow}

\newcommand{\wt}{\widetilde}

\newcommand{\pibar}{\overline{\Pi}}

\newcommand{\Xbar}{\overline{X}}
\newcommand{\Sbar}{\overline{S}}
\newcommand{\Gtau}{G_{\tau_u-}}
\newcommand{\R}{\mathbb{R}}

\newcommand{\topr}{\stackrel{\mathrm{P}}{\longrightarrow}}
\newcommand{\todr}{\stackrel{\mathrm{D}}{\longrightarrow}}
\newcommand{\eqdr}{\stackrel{\mathrm{D}}{=}}
\newcommand{\halmos}{$\sqcup\!\!\!\!\sqcap$}

\begin{document}
\title{Stability of the Exit Time for L\'evy Processes}

\authorone[Syracuse University]{Philip S. Griffin}
\addressone{Department of Mathematics, Syracuse University, Syracuse, NY, 13244-1150,USA}
\authortwo[Australian National University]
{Ross A. Maller}
\addresstwo{Centre for Financial Mathematics, and School of Finance,
Actuarial Studies, \& Applied Statistics, Australian National
University, Canberra, ACT, Australia}\footnote{Research partially
supported by ARC Grant DP1092502}

\begin{abstract}
This paper is concerned with the behaviour of a L\'{e}vy process when it
crosses over a  positive level, $u$, starting from 0, both as
$u$ becomes large and as $u$ becomes small.
Our main focus is on the time, $\tau_u$, it takes the process to
transit
above the level,
and in particular, on the  {\it stability} of  this passage time;
thus, essentially, whether or not $\tau_u$ behaves linearly as $u\dto
0$ or
$u\to\infty$.
We also consider conditional stability of $\tau_u$
when  the process drifts to $-\infty$, a.s.
This provides information
relevant to quantities associated with
the ruin of an insurance risk process, which we analyse under a
Cram\'er condition.
\end{abstract}

\bigskip

\bigskip

\setcounter{equation}{0} \section{Introduction}\label{s1}
For a random walk $S$ starting from 0 with a positive step length
distribution  and  having finite mean, the number of steps required to
first pass a positive level $u$, $\tau^S_u$, say, is, for large $u$,
asymptotic to a multiple of $u$, the constant of proportionality being
the reciprocal of the mean step length. More precisely, $\tau^S_u/u\to 1/c$ a.s. as $u\to \infty$,
 and, further, $E\tau^S_u/u\to 1/c$ as $u\to
\infty$, where $c\in(0,\infty)$ is the expected step length.
These express a kind of long term linearity of the passage
time, and provide useful intuition in applications. Together with
other ``fluctuation" quantities related to passage over a level, such
as the overshoot of the level, and various undershoots, etc., it
constitutes one of many well known properties of the renewal theory of
random walks. More generally, properties such as stability of the
passage time, etc., have been extended to random walks on the line.
(References and further discussion are given later.)

It is natural to consider carrying the discrete time results
over to a L\'evy process $(X_t)_{t\ge 0}$, and this has been done in the literature
for some of the fluctuation quantities, especially, see \cite{DM3}
for stability of the overshoot. Applications of this and related
kinds of result abound; we have in mind, in particular, applications
to the insurance risk process: see, e.g., recent results in
\cite{AK},
\cite{DK}, 
\cite{FZ},
\cite{kkm},  
and \cite{PM}.    
These authors have tended to concentrate on properties
of the overshoot and undershoots, with less attention paid to the
ruin time, $\tau_u$. But it could be argued that $\tau_u$ is the most important or at least the
most interesting variable, from a practical point of view.

Our aim in this paper is to set out in detail a comprehensive
listing of conditions for the stability of $\tau_u$, in the
L\'evy
setting. For ``large time'' stability, i.e., as $u\to\infty$, the
discrete time (random walk) results can be consulted to give useful
guidance for some of the L\'evy results; others are rather
straightforward to transfer, but others again are challenging. We
consider both stability in probability and almost sure (a.s.)
stability of $\tau_u/u$, as $u\to\infty$, when
$\limsup_{t\to\infty}X_t=\infty$ a.s., and
$\lim_{t\to\infty}X_t=-\infty$ a.s.

Even more interesting is the ``small time'' stability, i.e., as
$u\to 0$, of the passage time. Here there are of course no
corresponding random walks that can be used for guidance, but,
remarkably, small time results for L\'evy processes often parallel
large time results in certain ways. With this insight and some
further analysis we are able to give also a comprehensive analysis
of the small time stability of $\tau_u$. Some curious and unexpected
results occur (see, e.g., Remark \ref{2.3}). Such results may be thought of
as adding to our understanding of the local properties of L\'evy
processes.

The setting is as follows. Suppose that $X=\{X_{t}: t \geq 0 \}$,
$X_0=0$,  is a L\'{e}vy process defined on $(\Omega, {\cal F}, P)$,
with triplet $(\gamma, \sigma^2, \Pi_X)$, $\Pi_{X}$ being the
L\'{e}vy measure of $X$, $\gamma\in \mathbb{R}$,  and $\sigma^2\ge
0$. Thus the characteristic function of $X$ is given by the
L\'{e}vy-Khintchine representation, $E(e^{i\theta X_{t}}) = e^{t
\Psi_X(\theta)}$, where
\begin{equation}\label{lrep}
\Psi_X(\theta) =
 \rmi\theta \gamma - \sigma^2\theta^2/2+
\int_{\R}(e^{\rmi\theta x}-1-
\rmi\theta x \mathbf{1}_{\{|x|\le 1\}})\Pi_X(\rmd x),
\ {\rm for}\  \theta \in \mathbb{R}.
\end{equation}

Denote the  maximum process by
\[
\Xbar_{t} = \sup_{0\le s\leq t}X_s,
\]
and let
\[
G_t= \sup\{0\le s\le t: X_s= \Xbar_{s}\}
\]
be the time of the last maximum prior to time $t$.
Our focus will be on the
first passage time above level $u$, defined by
\[
 \tau_u= \inf \{t \geq0: X_t>u \},\ u>0.
\]
(We adopt the convention that the inf of the empty set is $+\infty$.)
Also important will be the
time of the last maximum before passage, $G_{\tau_u-}$, and the
position after  transit
above level $u$, $X_{\tau_u}$.
Throughout, we assume that  $\Pi_X$ is not identically zero
and that $X$ is not the negative of a subordinator
(in which case $\tau_u=\infty$ for all $u>0$).
By a  compound Poisson process we will mean a
L\'evy process with finite L\'evy measure, no Brownian component and
zero drift.

We need some further notation.
Let  $(L^{-1}_t,H_t)_{t \geq 0}$
denote the bivariate ascending
inverse local time--ladder height subordinator process of $X$.
The process $(L^{-1},H)$ is  defective when, and only when,
$\lim_{t\to \infty} X_{t} = -\infty$ a.s.
In that case, it
is obtained from  a nondefective process
$({\cal L}^{-1},{\cal H})$ by exponential killing with
rate $q > 0$, say.
When $(L^{-1},H)$ is nondefective the
killing is unnecessary and we set
$(L^{-1},H)=({\cal L}^{-1},{\cal H})$ and take $q=0$.
We denote the bivariate L\'{e}vy measure of
$({\cal L}^{-1}_t,{\cal H}_t)_{t\ge 0}$
by $\Pi_{{ L}^{-1},{ H}}(\cdot,\cdot)$,
and let $ \Pi_{L^{-1}}$ and $\Pi_H$  be the marginal
L\'evy measures of ${\cal L}^{-1}$ and ${\cal H}$.
The Laplace exponent $\kappa(a,b)$
of  $(L^{-1},H)$
will play an important role in our analysis.
It is defined by
\begin{equation} \label{kapdef}
e^{-\kappa(a,b)} = e^{-q}E e^{-a{\cal L}^{-1}_1 -b{\cal H}_1}
 \end{equation}
for values of $a,b\in \R$ for which the expectation  is finite.
We can write
\begin{eqnarray} \label{kapexp}
\kappa(a,b)
= q+\rmd_{L^{-1}}a+\rmd_Hb+\int_{t\ge 0}\int_{h\ge 0}
\left(1-e^{-at-bh}\right)
\Pi_{{ L}^{-1}, { H}}(\rmd t, \rmd h),
\end{eqnarray}
where $\rmd_{L^{-1}}\ge 0$ and $\rmd_H\ge 0$ are drift constants.
See, e.g.,
\cite{Bert}, \cite{doneystf}, \cite{kypbook}, \cite{sato},
for these relationships.

The following theorem connects the Laplace transform of the
fluctuation quantities with the bivariate Laplace exponent. It is an
extension of the ``second factorisation identity" (\cite{prog}, Eq.
(3.2)). A proof of Theorem \ref{thmct3}  is in \cite{GM2010}.

\begin{thm}[Laplace Transform Identity]\label{thmct3}\
Fix $\mu> 0$, $\rho\ge 0$, $\lambda\ge 0$, $\nu\ge 0$, $\theta\ge
0$.
If $\mu+\lambda\ne \rho$,
\begin{equation}\label{e3}
\int_{u\ge 0}e^{-\mu u}
E\left(e^{-\rho(X_{\tau_u}-u)-\lambda(u-\Xbar_{\tau_u-})
-\nu \Gtau-\theta(\tau_u-\Gtau)}; \tau_u<\infty\right)\rmd u
=\frac {\kappa(\theta,\mu+\lambda)-\kappa(\theta,\rho)}
{(\mu+\lambda-\rho)\kappa(\nu,\mu)}.
\end{equation}
 \end{thm}
In the present paper we apply these concepts to study the {\it
stability}
of the passage time, $\tau_u$,
by which we mean that $\tau_u/u$ has a finite and  positive
nonstochastic limit,
where the convergence may be as $u\to 0$ or $u\to\infty$, and
the convergence may be in probability, almost sure (a.s.), or in mean.
We will also consider, to a lesser extent,
 the position, $X_{\tau_u}$, of $X$ as it crosses
the boundary.
Some other results of interest, especially, that $\tau_u/u$ are
uniformly integrable as $u\to\infty$ if $X$ has a finite positive mean (see Lemma
\ref{Lailem}) are derived as by-products.

The results relating to stability of $\tau_u$ are given in Section \ref{s2}.
In contrast, in Section \ref{s3} we consider  large time conditional stability of $\tau_u$ when $P(\tau_u<\infty)\to 0$ as $u\to\infty$.
This is the usual
setup in the L\'evy insurance risk model, for which see, e.g.,
\cite{aa}, \cite{a02}, \cite{DK} and \cite{kkm}
for background and references. Section \ref{s4} contains some
concluding remarks and references.  All proofs are in Sections
\ref{s5}, \ref{s6}, and the Appendix.

\setcounter{equation}{0}
\section{Stability
}\label{s2}

\noindent
This section contains results relating to the stability of $\tau_u$ as $u\to L$ where $L=\infty$ or $L=0$. For stability  to make sense when $L=\infty$ we need, at a minimum, to assume that $P(\tau_u<\infty)\to 1$ as $u\to \infty$.  This is equivalent to $\limsup_{t\to\infty}X_t=+\infty$ a.s., in which case $\tau_u<\infty$ a.s. for all $u>0$ and $\tau_u\to \infty$ a.s. as $u\to\infty$.  The natural analogue of this condition when $L=0$ is that $P(\tau_u<\infty)\to 1$ and $\tau_u\to 0$ a.s. as $u\dto 0$.  This is equivalent to $0$ being regular for $(0,\infty)$;
see \cite{bert2} for an analytic equivalence. Thus {\it the overriding assumptions throughout this section are: $\limsup_{t\to\infty}X_t=+\infty$ a.s. when $L=\infty$, and $0$ is regular for $(0,\infty)$ when $L=0$.}

Let $\pibar_X$ and $\pibar_X^\pm$ denote the tails of $\Pi_X$, thus
\begin{equation}
\pibar_X^+(x)=\Pi_X\{(x,\infty)\},\ \pibar_X^-(x)=\Pi_X\{(-\infty, -
x)\},\ {\rm and}\
 \pibar_X(x)= \pibar_X^+(x)+\pibar_X^-(x),
\end{equation}
for $x>0$, and define a kind of truncated mean
 \begin{eqnarray}\label{Adef}
A(x)&:=& \gamma +\pibar_X^+(1)- \pibar_X^-(1)
+\int_1^x \left(\pibar_X^+(y)- \pibar_X^-(y)\right)\rmd y
\nonumber\\
&=&
\gamma + x\left(\pibar_X^+(x)- \pibar_X^-(x)\right)
+\int_{1<|y|\le x} y \Pi_X(\rmd y),\ x>0.
 \end{eqnarray}

The first theorem concerns the stability in probability of $\tau_u$.
Consider first stability for large times, as $u\to\infty$, i.e., the
property $\tau_u/u\topr 1/c$ as $u\to\infty$, for some $c\in
(0,\infty)$. This is equivalent to the relative stability in
probability of the process $X$ itself, i.e., to $X_t/t\topr c$ as
$t\to\infty$. We prove it via an equivalence of the stability of
$\tau_u$ with that of $\Xbar$, namely, $\Xbar_t/t\topr c$ as
$t\to\infty$, a trivial relationship.  We  then show that the latter
holds iff $X$ is relatively stable, which is not entirely obvious,
but follows from similar (large time) random walk working of
\cite{KM99}, where the stability of the passage time of a random
walk above a constant level is considered for general norming
sequences. The stability of $\tau_u$ is connected to the bivariate
Laplace exponent in \eqref{rt2}, which is a new relationship,
derived via Theorem \ref{thmct3}, and the list of equivalences for
this case is completed by that of \eqref{rt5} and \eqref{rt6}, which
is in Theorem 3.1 of  \cite{DM2}.

This list, for the case $u\to\infty$, $c\in (0,\infty)$, then sets the
pattern we work from for the case $u\dto 0$, $c\in (0,\infty)$, and later results.
Theorem \ref{taurs} also considers the cases $c=0$ and $c=\infty$ for
completeness, though these strictly speaking do not give rise to
stability conditions.

\begin{thm}[Stability in Probability of the Exit Time]\label{taurs}
(a)\
Fix a constant  $c\in(0,\infty)$ and
let $L=0$ or $\infty$ ($1/L=\infty$ or $0$).
Then the  following are equivalent:
\begin{equation}\label{rt1}
\frac{\tau_u}{u} \topr \frac{1}{c},\ {\rm as}\ u\to L;
\end{equation}
\begin{equation}\label{rt2}
\lim_{x\to 1/L}\frac{\kappa(x,0)}{\kappa(x,\xi x)}=\frac{1}{1+\xi c},
 \  {\rm for\ each}\ \xi>0;
\end{equation}
\begin{equation}\label{rt41}
\frac{\Xbar_t}{t} \topr c,\ {\rm as}\ t\to L;
\end{equation}
\begin{equation}\label{rt5}
\frac{X_t}{t} \topr c,\ {\rm as}\ t\to L.
\end{equation}

In the case $L=\infty$, \eqref{rt1}--\eqref{rt5} are equivalent to
 \begin{equation}\label{rt6}
x\pibar_X(x) \to 0\ {\rm and}\ A(x)\to c, \  {\rm as}\ x\to\infty.
 \end{equation}

In the case $L=0$, \eqref{rt1}--\eqref{rt5} are equivalent to
 \begin{equation}\label{rt46}
\sigma^2=0, \
x\pibar_X(x) \to 0\ {\rm and}\ A(x)\to c, \  {\rm as}\ x\dto 0.
 \end{equation}

 \noindent
(b)\
Suppose $c=0$. If $L=\infty$ then
\eqref{rt1}--\eqref{rt6} remain equivalent.  If $L=0$ then \eqref{rt1}--\eqref{rt41} remain equivalent, as do \eqref{rt5} and \eqref{rt46}.  However while  \eqref{rt5} implies \eqref{rt1}--\eqref{rt41}, the converse does not hold.

\noindent
(c)\
Suppose $c=\infty$. Then
\eqref{rt1}--\eqref{rt41}
remain equivalent for $L=0$ or $\infty$ in the following sense:
 \begin{equation}\label{rt0}
\frac{\tau_u}{u} \topr 0,\ {\rm as}\ u\to L,
\ {\rm iff}\
\lim_{x\to 1/L}\frac{\kappa(x,0)}{\kappa(x,\xi x)}=0,  \  {\rm for\
each}\ \xi>0,
\ {\rm iff}\
\frac{\Xbar_t}{t} \topr \infty,\ {\rm as}\ t\to L.
\end{equation}
Again, while
\eqref{rt5} implies \eqref{rt1}--\eqref{rt41}, it is
not equivalent
in either case,  $L=0$ or $\infty$.
\end{thm}

\begin{rem}
As mentioned above, the equivalence of  \eqref{rt5} and \eqref{rt6} when $L=\infty$
is in \cite{DM2}, while the  equivalence of \eqref{rt5} and \eqref{rt46} when $L=0$ is in Theorem 2.1 of \cite{DM2}.  Both of these results hold for all $c\in (-\infty,\infty)$.
We include them in the statement of Theorem \ref{taurs} for completeness and for the convenience of the reader.
\end{rem}

The next theorem concerns the almost sure stability of $\tau_u$. We
follow the pattern set by Theorem \ref{taurs}. The connection with the
bivariate Laplace exponent is transmuted in this case to requiring
finite first moments  of the ladder processes $H$ and $L^{-1}$.
Almost sure stability for large times requires a finite positive mean for $X$ (for
large times), and bounded variation with positive drift of $X$ (for small times). Recall that when $X$ is of bounded variation, we may write
 the L\'evy-Khintchine exponent in the form
\be\label{lkbv}
\Psi(\theta)=\rmi\theta\rmd_X+
\int_{\R}(e^{\rmi\theta x}-1)\Pi_X(\rmd x),
\ee
where
$
\rmd_X=\gamma-\int x \mathbf{1}_{\{|x|\le 1\}}\Pi_X(\rmd x)
$
is called the drift of $X$.

\begin{thm}[Almost Sure Stability of the Exit Time]\label{taursas}

\noindent
(a)\ Fix $c\in[0,\infty)$. \
(i)\
We have  $\displaystyle{\frac{\tau_u}{u}} \to \frac{1}{c}$, a.s.
as $u\to \infty$ iff $E|X_1|<\infty$ and  $EX_1=c\ge 0$.

(ii)\
We have  $\displaystyle{\frac{\tau_u}{u}} \to \frac{1}{c}$, a.s.,
as $u\to 0$ iff
$X$ is of bounded variation with drift  $\rmd_X=c\ge 0$.

\noindent
(b)\ Fix $c\in(0,\infty)$.
Then (i) holds iff $EH_1<\infty$ and $EL_1^{-1}<\infty$,
in which case
 $c=EH_1/EL_1^{-1}$, while (ii) holds iff $\sigma^2=0$,
${\rm d}_{L^{-1}}>0$   and  ${\rm d}_{H}>0$,
in which case $c= {\rm d}_{H}/{\rm d}_{L^{-1}}$.
\end{thm}

\begin{rem}
Note that, under \eqref{rt6} and \eqref{rt46} respectively,
\begin{equation}\label{Aint}
\lim_{x\to \infty} A(x)=\gamma+\int_{|y|>1} y \Pi_X(\rmd y)
\quad {\rm and}\quad
\lim_{x\dto 0} A(x)=\gamma-\int_{0<|y|\le 1} y \Pi_X(\rmd y).
\end{equation}
Here existence of the limits is equivalent to conditional convergence of the integrals.  Under the conditions of parts
(a)(i) and (a)(ii) of Theorem \ref{taursas}, these integrals converge absolutely and the limits are then given by $EX_1$
and $\rmd_X$ respectively,
thus confirming that  the expressions for $c$ in
Theorems \ref{taurs} and \ref{taursas} agree.
The difference between \eqref{rt6} and \eqref{rt46} and  (i) and (ii) of Theorem
\ref{taursas}  is essentially whether the integrals in \eqref{Aint}
converge conditionally or absolutely.

\end{rem}

In the next theorem we examine the convergence of $E\tau_u/u$ as
$u\to\infty$ and as $u\dto 0$.
Recall that  $E\tau_u<\infty$ for some, hence all, $u\ge 0$,
iff $X$ drifts to $+\infty$ a.s., iff $EL_1^{-1}<\infty$
(see, e.g.,  Theorem 1 of  \cite{DM4}).

\begin{thm}[Stability of the Expected Exit Time]\label{Etau}
(a)\
Fix $c\in(0,\infty)$. Then

(i)\ $E\tau_u<\infty$ for each $u>0$ and
$\displaystyle{\lim_{u\to\infty}\frac{E\tau_u}{u}}=\frac{1}{c}$
iff $0<EX_1\le E|X_1|<\infty$. In this situation,
$EH_1<\infty$,  $EL_1^{-1}<\infty$, and $c=EX_1= EH_1/EL_1^{-1}$.

(ii)\  $E\tau_u<\infty$ for each $u>0$ and
$\displaystyle{\lim_{u\dto 0}\frac{E\tau_u}{u}} =\frac{1}{c}$
iff
$EL_1^{-1}<\infty$ and $\rmd_H>0$,
and then $c=\rmd_H/EL_1^{-1}$.

\noindent
(b)\ (The case $c=0$)
(i)\
In Part (i) of  the theorem, the case $c=0$ cannot arise;
 when $E\tau_u<\infty$ for each $u>0$,
 $\lim_{u\to\infty}E\tau_u/u$ exists and is in $[0,\infty)$.

(ii)\
We have $E\tau_u<\infty$ for each $u>0$, and
$\lim_{u\dto 0}E\tau_u/u=\infty$ iff  $EL_1^{-1}<\infty$ and
$\rmd_H=0$.
\noindent

\noindent
(c)\ (The case $c=\infty$)
(i)\
We have $E\tau_u<\infty$ for each $u>0$ and
$\lim_{u\to \infty}E\tau_u/u=0$ iff   $EL_1^{-1}<\infty$ and
$EH_1=\infty$.

(ii)\
In Part (ii) of  the theorem, the case $c=\infty$ cannot arise;
 when $E\tau_u<\infty$ for each $u>0$, we always have $\liminf_{u\dto
0}E\tau_u/u>0$.
\end{thm}

\begin{rem}\label{2.3}
(i)\
It is curious that the formula for $c$ in Part (a)(ii) of
Theorem \ref{Etau} doesn't agree with the versions
in   Theorem \ref{taurs} or   Theorem \ref{taursas}.
We give an example to illustrate how the difference can arise.
Let
\ben
X_t=at-N_t
\een
where $N_t$ is a rate one Poisson process and $a>1$.
Thus $\lim_{t\to\infty}X_t=\infty$ a.s.
Since $\tau_u=ua^{-1}$ for sufficiently small $u$,
it trivially follows that
\ben
\lim_{u\dto 0}{\frac{\tau_u}{u}} =\frac{1}{a}\ {\rm a.s.}
\een
We claim that
\be\label{etucp}
{\frac{E\tau_u}{u}} \to \frac{1+E\tau_1}a.
\ee
This is because, if $\xi$ is the time of the first jump of $N$, then
\ben\ba
E\tau_u&=E(\tau_u;N_{ua^{-1}}=0)+E(\tau_u;N_{ua^{-1}}\ge 1)\\
&={ua^{-1}}e^{-{ua^{-1}}} + \int_0^{ua^{-1}}E(\tau_u|\xi =t)P(\xi\in \rmd t)\\
&={ua^{-1}}e^{-{ua^{-1}}} + \int_0^{ua^{-1}} (t+E\tau_{1-at+u})e^{-t} \rmd t.
\ea\een
Since $\tau_{1+x}\overset{P}\to\tau_{1}$ as $x\dto 0$,
\eqref{etucp} now follows after dividing by $u$ and taking the limit.

We now check that this agrees with Part (ii) of Theorem \ref{Etau}.
For the normalisation of $L$, the local time at the maximum,  we take
\ben
L_t=\int_0^t I(X_s = \overline{X}_s)\ \rmd s.
\een
Then the  ladder height process
is linear drift, ${H}_{t}=at$.  Hence $\rmd_H=a$.  By construction
\ben
L_t^{-1}=t+\sum_1^{N_t} R_i
\een
where again $N_t$ is a rate one Poisson process and $R_i$ are iid
 random variables independent of $N$, with distribution the same as that of
$\tau_1$.  Hence  $\rmd_{L^{-1}}=1$ and
\ben
EL_1^{-1}=1+E\tau_1,
\een
so the  $c$ in Part (ii) of Theorem \ref{Etau} is
$\rmd_H/EL_1^{-1}=a/(1+E\tau_1)$, giving
agreement with \eqref{etucp}.
On the other hand,  the  $c$ in Part (ii) of Theorem \ref{taursas} is
$\rmd_X=\rmd_H/\rmd_{L^{-1}}=a$.
 \end{rem}

We now turn to stability of the time of the last maximum before ruin. As may be expected, this is a more difficult object to study than $\tau_u$.  We consider the three modes of convergence investigated in Theorems
\ref{taurs}--\ref{Etau}.

\begin{thm}[Stability of the Last Maximum before Ruin]
\label{Grs}
Let $L=0$ or $\infty$ ($1/L=\infty$ or $0$).

\noindent(a)\ Fix $c\in (0,\infty)$.
 We have
\begin{equation}\label{rt4}
\frac{G_{\tau_u-}}{u} \topr \frac{1}{c},\ {\rm as}\ u\to L,
\end{equation}
if and only if
\begin{equation}\label{rt3}
\lim_{x\to 1/L}\frac{\kappa(0,x)}{\kappa(\xi x,x)}=\frac{c}{c+\xi},
\  {\rm for\ each}\ \xi>0.
\end{equation}
(b) Fix $c\in [0,\infty)$.
(i)\
$\displaystyle{\frac{G_{\tau_u-}}{u}} \to \frac{1}{c}$, a.s.
as $u\to \infty$ iff $E|X_1|<\infty$ and  $EX_1=c\ge 0$.

(ii)\
$\displaystyle{\frac{G_{\tau_u-}}{u}} \to \frac{1}{c}$, a.s.,
as $u\to 0$ iff
$X$ is of bounded variation with drift  $\rmd_X=c\ge 0$.

\noindent(c)\
 If $0<EX_1<\infty$ then
$\displaystyle{\lim_{u\to\infty}\frac{E G_{\tau_u-}}{u}}=\frac{1}{EX_1}$.
%
\end{thm}

\begin{rem}
It's not clear how the conditions of Theorem \ref{taurs}
relate to the stability in probability of $G_{\tau_u-}$.
We can show that \eqref{rt1}--\eqref{rt6} imply  \eqref{rt4} and
\eqref{rt3} but it's not clear whether or not the converse holds.  For almost sure convergence  the results for
$G_{\tau_u-}$ parallel those for ${\tau_u}$. For convergence in mean the situation remains largely unresolved.
\end{rem}

The final result, Theorem \ref{jet}, belongs in the present section since it
holds in  a case when $\lim_{t\to\infty}X_t=+\infty$ a.s.,
but we apply it in  the next section, in
the case when $\lim_{t\to\infty}X_t=-\infty$ a.s.,
to obtain results in the L\'evy insurance risk model.

\begin{thm}[Convergence of Expected Exit Times with Overshoot]
\label{jet}\
\newline
Assume $0<EX_1\le E|X_1|<\infty$, and  that  $X$ is not compound Poisson, or is
compound Poisson with a nonlattice jump distribution.
Then for all $\rho>0$
\begin{equation}\label{to1}
\lim_{u\to\infty}E\left(\frac{G_{\tau_u-}}{u}e^{-\rho(X_{\tau_u}-u)}
\right)
=\lim_{u\to\infty}E\left(\frac{\tau_u}{u}e^{-\rho(X_{\tau_u}-
u)}\right)
=\frac{1}{EX_1} E e^{-\rho Y},
\end{equation}
where $Y$ is the limiting distribution of the overshoot $X_{\tau_u}-u$.  $Y$ has density
$\pibar_H(h)\rmd h/EH_1$ on $(0,\infty)$, and mass $\rmd_H/EH_1$ at 0.
\end{thm}

\setcounter{equation}{0}
 \section{Stability In the Insurance Risk Model}\label{s3}

The aim of this section is to illustrate that stability questions are also of interest when $X_t\to -\infty$ a.s. as $t\to\infty$.
We phrase the discussion in terms  of an insurance risk model. In
this case $X$ represents the excess in claims over premium  of an
insurance company. The classical model in this context is the
Cram\'er-Lundberg model in which  $X$ is the sum of a compound
Poisson process with positive jumps, representing claims, and a
negative drift, representing premium inflow.  The results in the
present section will be given for a general L\'evy insurance risk
model where  no such restrictions are placed on $X$.

The over-riding assumption throughout  this section  is the {\em
Cram\'er condition}, namely, that
\begin{equation}\label{cramer}
Ee^{\nu_0X_1}=1, \ {\rm  for\ some}\  \nu_0>0.
\end{equation}
It's well known that, under  \eqref{cramer}, $EX_1$ is well defined,
with $EX_1^+<\infty$, $EX_1^-\in(0,\infty]$, and $EX_1\in
[-\infty,0)$, and  so $\lim_{t \to \infty} X_t=-\infty$ a.s.
Further,
$E(X_1e^{\nu X_1})$ is finite and positive for all $\nu$ in a left neighbourhood
of $\nu_0$, and
\be\label{mu*>0}
\mu^*:= E(X_1e^{\nu_0 X_1})>0\ \   (\rm{possibly }\ \mu^*=+\infty).
\ee


Since   $\lim_{t \to \infty} X_t=-\infty$ a.s., we are in the
situation that $P(\tau_u<\infty)<1$ for all $u>0$,
and $\lim_{u\to\infty}P(\tau_u<\infty)=0$.
In an insurance risk context, we are interested in forecasting
the ruin time $\tau_u$ in a worst case scenario,
i.e., conditional on $\tau_u<\infty$ (``ruin occurs").
Asymptotic properties of $\tau_u$ and associated variables,
conditional on $\tau_u<\infty$, often provide surprisingly
good approximations of corresponding
finite level distributions; cf, e.g., \cite{grandell}.
In the present context we look at the
stability of $\tau_u$, $G_{\tau_u-}$ and $X_{\tau_u}$, showing they are asymptotically
linear under mild conditions.

We need some more infrastructure. Let  $(X^*_t)_{t\ge 0}$ denote the
Esscher transform of $X$ defined by
 \begin{equation}\label{X*X}
P\left( (X^*_s, 0\le s\le t)\in B, X^*_t\in \rmd x\right)
=e^{\nu_0 x}
P\left( (X_s, 0\le s\le t)\in B, X_t\in \rmd x\right),
 \end{equation}
for any Borel subset $B$ of $\R^{[0,t]}$.
Equivalently $X^*$  may be introduced by means of exponential tilting;  that is,
define a new probability $P^*$, given on ${\cal F}_t$ by
\be
\frac{dP^*}{dP}=e^{\nu_0 X_t}.
\ee
Then $X$ under $P^*$ has the same distribution as $X^*$ under $P$.
It easily follows that
\begin{equation*}
Ef(X_t^*)=E^*f(X_t)=E(f(X_t)e^{\nu_0 X_t}),
 \end{equation*}
for any Borel function $f$ for which the expectations are finite.
$X^*$ is itself a L\'evy process with exponent  $\Psi(\theta
-\rmi\nu_0)$ and $E^*X_1=\mu^*$. Since $\mu^*>0$ by  \eqref{mu*>0},
$X^*_t$ drifts to $+\infty$ a.s., and hence $(H^*_t)_{t\ge 0}$, the
increasing ladder height process associated with $(X^*_t)_{t\ge 0}$,
is proper.

Our setup is that of Bertoin and Doney \cite{BD94}.
The main result in \cite{BD94},
which we give in  the form proved in Theorem 7.6 of \cite{kypbook},
see also  Section XIII.5 of \cite{a02}, is:
{\it Suppose \eqref{cramer} holds and
the support of $\Pi_X$ is non-lattice in the case that
$X$ is compound Poisson. Then
 \begin{equation}\label{BDmain}
\lim_{u\to\infty}e^{\nu_0 u}P(\tau_u<\infty)=C\in[0,\infty),
 \end{equation}
where $C:=E^*e^{-\nu_0Y}>0$ if and only if $\mu^*<\infty$.
Here  $Y$ is the limiting distribution of the overshoot $X_{\tau_u}-u$ under $P^*$.
}

To state the stability result for the general L\'evy insurance risk
model under \eqref{cramer}, introduce the probability  measure
$P^{(u)}(\ \cdot\ )=P(\ \cdot\ |\tau_u<\infty)$, and denote
convergence in probability conditional on $\tau_u<\infty$ by
$\topru$.

\begin{thm}\label{thm1}
Assume \eqref{cramer} holds and $\mu^*<\infty$ (so that  $0<\mu^*<\infty$).
Then,  as $u\to\infty$,
\begin{equation}\label{Astau}
\frac{X_{\tau_u}}{u}\topru 1,\ \frac{G_{\tau_u-}}{u}\topru \frac{1}{\mu^*}
\quad {\rm and}\quad
\frac{\tau_u}{u}\topru \frac{1}{\mu^*}.
\end{equation}
Assume in addition that the support of $\Pi_X$ is non-lattice in the case that
$X$ is compound Poisson.
Then
 \begin{equation}\label{Eu0}
\lim_{u\to\infty}\frac{E^{(u)}X_{\tau_u}}{u}=1,\ \lim_{u\to\infty}
\frac{E^{(u)}G_{\tau_u-}}{u}=\frac{1}{\mu^*}
 \quad {\rm and}\quad
\lim_{u\to\infty}
\frac{E^{(u)}\tau_u}{u}=\frac{1}{\mu^*}.
\end{equation}
\end{thm}

Parts of our Theorem \ref{thm1} are well known for the
Cram\'er-Lundberg model, and their extension to the general L\'evy
insurance risk model is straightforward.  Others
appear to be new.

\setcounter{equation}{0} \section {Concluding Remarks}\label{s4}
There is of course a very large literature on (large time) renewal
theorems for random walks, and, more recently, some similar results
have been proved for L\'evy processes. Regarding the ruin time, most
results so far concern the infinite horizon ruin probability,
$P(\tau_u<\infty)$, or, equivalently, the distribution of the
overall maximum of the random walk or L\'evy process, and we do not
attempt to summarise them here (other than the references mentioned
in Sections \ref{s1}--\ref{s3}). A web search turns up many such
papers and books.

The finite horizon ruin probability, $P(\tau_u<T)$, is less studied,
but important results are obtained in, e.g., \cite{A84}, \cite{AK},
\cite{HDCW}, \cite{Avrametal}, \cite{haotang},
 \cite{Brave} (see also their references), and especially, in the
 insurance/actuarial literature (usually from a more applied point of view). These results of
course give information on the long run distribution of the ruin
time, conditional on ruin occurring. A more recent result along
these lines is in \cite{GM2}, assuming, like  \cite{Brave}  and
\cite{haotang}, convolution equivalent conditions on the tails of
the process or its L\'evy measure. These authors are interested in
the asymptotic distribution of $\tau_u$, rather than in its
stability {\em per se}; as mentioned earlier, results on stability
such as we give are more akin to classical (large time) renewal
theory than to these, and small time versions, which make sense for
L\'evy processes but not for random walks, have previously been
neglected, in the main.

We turn now to the proofs.

\setcounter{equation}{0} \section{Proofs for Section
\ref{s2}}\label{s5}
We assume throughout this section that $\limsup_{t\to\infty}X_t=+\infty$ a.s. when $L=\infty$, and
$0$ is regular for $(0,\infty)$ when $L=0$.  In the former case, $\tau_u<\infty$ a.s. for all $u>0$, while $P(\tau_u<\infty)\to 1$ as $u\to 0$ in the latter case.

\bigskip \noindent {\bf Proof of Theorem \ref{taurs}}.\
(a) and (b).\
Let $c\in [0,\infty)$ until further notice, with the obvious
interpretations when $c=0$.
Assume \eqref{rt1} holds with $u\to L$.
From \eqref{e3} we have, for $u>0$, $\theta>0$,
\begin{equation}\label{tauform}
\mu\int_{u\ge 0}  e^{-\mu u}
E\left(e^{-\theta\tau_u};\tau_u<\infty\right)\rmd u
=1-\frac{\kappa(\theta,0)}{\kappa(\theta,\mu)}.
\end{equation}
Take $y>0$ and replace $\mu$ by $\mu/y$, $u$ by $uy$, and $\theta$ by $\theta/y$ in this to get
 \begin{equation}\label{rs1}
\mu\int_{u\ge 0}  e^{-\mu u}
E\left(e^{-\theta\tau_{uy}/y};\tau_{uy}<\infty\right)\rmd u
=1-\frac{\kappa(\theta/y,0)}{\kappa(\theta/y,\mu/y)}.
\end{equation}
By hypothesis, $\tau_{uy}/y\topr u/c$ as $y\to L$ for each $u>0$,
so letting  $y\to L$ in \eqref{rs1} gives, by dominated convergence,
 \begin{equation*}
\lim_{y\to L}
\frac{\kappa(\theta/y,0)}{\kappa(\theta/y,\mu/y)}
=1-  \mu\int_{u\ge 0}  e^{-\mu u} e^{-\theta u/c} \rmd u
=\frac{\theta}{\theta+\mu c}.
\end{equation*}
Replacing $y$ by $x=\theta/y\to 1/L$ and $\mu/\theta$ by $\xi$ gives
\eqref{rt2} with $x\to 1/L$.

Conversely, assume \eqref{rt2} with $x\to 1/L$.
Then from \eqref{rs1} we see that, for $\theta>0$,
\begin{equation*}
 \lim_{y\to L}
\mu\int_{u\ge 0}  e^{-\mu u}
E\left(e^{-\theta\tau_{uy}/y};\tau_{uy}<\infty\right)\rmd u
= \frac{\mu c}{\theta+\mu c}.
\end{equation*}
For each $y>0$, $\theta>0$, the function
$f_y(u,\theta):=  E\left(e^{-\theta\tau_{uy}/y};\tau_{uy}<\infty\right)$
is monotone decreasing in $u$ and bounded by 1.
Given any sequence $y_k\to L$ we can by Helly's theorem
find a subsequence  $\wt y_k\to L$, possibly depending on $\theta$
but not on $u$, such that
$f_{\wt y_k}(u,\theta)\to \wt f(u,\theta)$ for some function
$\wt f(u,\theta)\in [0,1]$.
Then by dominated convergence we have
 \begin{equation*}
\mu\int_{u\ge 0}  e^{-\mu u}\wt f(u,\theta)\rmd u
 = \frac{\mu c}{\theta+\mu c}
=\mu\int_{u\ge 0}   e^{-\mu u} e^{-\theta u/c} \rmd u,
\end{equation*}
and from the uniqueness of Laplace transforms we deduce that
$\wt f(u,\theta)=e^{-\theta u/c}$, not dependent on the choice
of subsequence.  Hence (taking $u=1$ now)
\[
\lim_{y\to L} E\left(e^{-\theta\tau_{y}/y};\tau_{y}<\infty\right)= e^{-\theta/c},
\ \theta>0,
\]
proving  \eqref{rt1} with $u\to L$.

Since
\begin{equation}\label{invr}
\{\tau_u>t\}  \subseteq \{\Xbar_t\le u\}  \subseteq \{\tau_u\ge t\},
\ t>0, u>0,
\end{equation}
we easily see that \eqref{rt1} is equivalent to
\eqref{rt41} in either case,  $L=\infty$ or $L=0$, for $c\ge 0$.

Next we show \eqref{rt5} implies \eqref{rt41}.  First consider  the case $L=\infty$.
By Theorem 3.1 of \cite{DM2},  \eqref{rt5} with $t\to\infty$
is equivalent to
\eqref{rt6}, and the first relation in \eqref{rt6} implies
\begin{equation}\label{rso}
\lim_{x\to \infty}\frac{V(x)}{x} =0,
 \end{equation}
where
\[
V(x):=\sigma^2+\int_{|y|\le x}y^2 \Pi_X(\rmd y).
\]
To deduce \eqref{rt41} from \eqref{rt6} and \eqref{rso} in the case
$L=\infty$,
decompose $X$ into small and large jump components as in Lemma 6.1 of
\cite{DM2} to get
 \begin{equation*}
X_s=s\nu(t)+\sigma B_s+X_s^{(1)}+X_s^{(2)},\ 0\le s\le t,
  \end{equation*}
where
\[
\nu(x):=\gamma  +\int_{1<|y|\le x} y \Pi_X(\rmd y),\ x>0,
\]
$\sigma B_s+X_s^{(1)}$ is a mean 0 martingale with jumps bounded in
modulus by $t$ and all moments finite,
and
\begin{equation*}
X_s^{(2)}=\sum_{0<r\le s}\Delta X_r {\bf 1}_{\{|\Delta X_r|>t\}}.
\end{equation*}
By Doob's inequality, for $\veps>0$,
\begin{equation*}
P\left(\sup_{0\le s\le t}|\sigma B_s+X_s^{(1)}|>\veps t\right)
\le \frac{1}{(\veps t)^2}E(\sigma B_t+X_t^{(1)})^2
=\frac{1}{\veps^2t}V(t),
\end{equation*}
and this tends to 0 as $t\to\infty$ by \eqref{rso}.
Also
\begin{equation*}
P\left(\sup_{0\le s\le t}|X_s^{(2)}|=0\right)
\geq P\left(\text{no jumps with }|\Delta X_s|>t\text{ occur
by time }t\right) =\exp (-t\pibar_X(t))\rightarrow 1,
\end{equation*}
as $t\to\infty$,
while
\[
\nu(t)=A(t)-t\pibar^+(t)+t\pibar^-(t) \to c,
\]
both by \eqref{rt6}. Thus we have
$P(\sup_{0\le s\le t}|X_s-cs|>\varepsilon t)\to 0$  as  $t\to\infty$,
which implies
\eqref{rt41} with $L=\infty$ and $c\ge 0$.

Next we deal with the implication
\eqref{rt5} implies \eqref{rt41} in the case  $L=0$.
Note that,
by Theorem 2.1 of \cite{DM2},  \eqref{rt5} with $t\dto 0$
is equivalent to \eqref{rt46},
and the first and second relations in \eqref{rt46} imply
\begin{equation}\label{rso2}
\frac{V(x)}{x}
=\frac{2}{x}\int_0^xy\pibar_X(y)\rmd y +x\pibar_X(x)
=\frac{2}{x}\int_0^xo(1)\rmd y +o(1)
=o(1),\ {\rm as}\ x\downarrow 0.
 \end{equation}
This takes the place of \eqref{rso} in the $L=0$ case, and the
rest of the proof that \eqref{rt5} implies \eqref{rt41}
with $L=0$ is virtually the same as for the case with $t\to\infty$.

We have left to show that \eqref{rt41} implies
 \eqref{rt5} except when $c=L=0$.
This is obvious if $c=0$ and $L=\infty$ so suppose $c>0$.
Note that, for $t>0$,
 \begin{eqnarray*}
\Xbar_{2t}&=& \Xbar_t  \vee \left(\sup_{t<s\le 2t}X_s\right)
\nonumber\\ &=&
\Xbar_t  \vee \left(X_t+\sup_{t<s\le 2t}(X_s-X_t)\right)
\nonumber\\ &=&
\Xbar_t  \vee \left(X_t+\Xbar_t'\right),
 \end{eqnarray*}
where $\Xbar_t'$ is an independent copy of $\Xbar_t$.
Consequently, for $\veps\in(0,c/3)$, as $t\to L$,
\begin{eqnarray*}
o(1)&=&P\left(\Xbar_{2t}\le (c-\veps)2t\right)
=P\left(\Xbar_t  \vee \left(X_t+\Xbar_t'\right)\le (c-\veps)2t\right)
 \nonumber\\ &\ge&
P\left(\Xbar_t \le (c-\veps)2t, X_t+(c+\veps)t\le (c-\veps)2t,
|\Xbar_t'-ct|\le \veps t\right)
 \nonumber\\ &=&
P\left(\frac{\Xbar_t}{t}\le 2(c-\veps),\frac{X_t}{t}\le c-3\veps\right)
P\left(\left|\frac{\Xbar_t}{t}-c\right|\le \veps\right)
 \nonumber\\ &\ge&
\left( P\left(\frac{X_t}{t}\le c-3\veps\right)
-  P\left(\frac{\Xbar_t}{t}>2(c-\veps)\right)\right)
\left( 1+o(1)\right)
  \nonumber\\ &=&
\left( P\left(\frac{X_t}{t}\le c-3\veps\right)-o(1)\right)
\left( 1+o(1)\right).
 \end{eqnarray*}
This shows that $P(X_t>(c-3\veps)t)\to 1$, and since also
\[
P\left(\frac{X_t}{t}\le c+\veps\right) \ge
P\left(\frac{\Xbar_t}{t}\le c+\veps\right) \to 1,
\]
we have $X_t/t\topr c$.
Hence \eqref{rt41} implies \eqref{rt5} for $L=0$ or $L=\infty$.

\medskip
(c)\
The equivalences in \eqref{rt0} follow by the same methods as used in Part (a),
and clearly \eqref{rt5} in case $c=\infty$ implies
\eqref{rt41} in case $c=\infty$.

All that remains is to give counterexamples showing  \eqref{rt41} does not imply \eqref{rt5} when $c=L=0$ or when $c=\infty$ and $L=0$ or $\infty$.

\begin{lem}\label{cex}
There is a L\'evy process for which 0 is regular for $(0,\infty)$  and with
\begin{equation}\label{ce1}
\frac{\Xbar_t}{t} \topr 0,\ {\rm as}\ t\dto 0,
\quad  {\rm but\  with}\quad
\frac{X_t}{t} \topr -1,\ {\rm as}\ t\dto 0.
\end{equation}
There is also a L\'evy process with
\begin{equation}\label{ce}
\frac{\Xbar_t}{t} \topr \infty,\ {\rm as}\ t\to L,
\quad  {\rm but\  with}\quad
\frac{X_t}{t} \topr -\infty,\ {\rm as}\ t\to L,
\end{equation}
for $L=0$ or $L=\infty$.
\end{lem}

\bigskip \noindent {\bf Proof of Lemma \ref{cex}}.\
This is given in the Appendix.
\hfill\halmos

With  Lemma \ref{cex} we complete the proof of Theorem \ref{taurs}.
\hfill\halmos

\bigskip \noindent {\bf Proof of Theorem \ref{taursas}}.\ By a simple pathwise argument using \eqref{invr}, it easily follows that
for $L=\infty$ or $L=0$,  and $c\in[0,\infty]$, we have $\displaystyle \lim_{u\to L}\frac{\tau_u}{u} = \frac{1}{c}\ {\rm  a.s.}$ if and only if
\begin{equation}\label{tauXbar}
\lim_{t\to L}\frac{\Xbar_t}{t} = c\ {\rm  a.s.}
 \end{equation}

(a)(i) If $E|X_1|<\infty$ and $c=EX_1\ge 0$ then \eqref{tauXbar} with $L=\infty$ holds by the strong law.
Conversely, by Theorem 15 of \cite{doneystf},  \eqref{tauXbar} implies that at least one of $EX_1^+$ or $EX_1^-$ is finite or else $X_t\to -\infty$.
Since $\limsup_{t\to\infty}X_t=\infty$, the latter possibility is ruled out and so is the possibility that $EX_1^-=\infty$.   Since $c\in[0,\infty)$, the strong law and \eqref{tauXbar} then force $E|X_1|<\infty$ and $c=EX_1\ge 0$.

(a)(ii) If $X$ is of bounded variation with $\rmd_X\ge 0$ then $X_t/t\to \rmd_X$ a.s. as $t\dto 0$ by Theorem 39 of \cite{doneystf}.  Hence \eqref{tauXbar} holds with $c=\rmd_X$.  Conversely, by the same result, \eqref{tauXbar}  implies  $X$ is of bounded variation and necessarily $\rmd_X\ge 0$ since $0$ is regular for $(0,\infty)$.  It then follows that $c=\rmd_X\ge 0$.

(b)(i)
When  $0<EX_1\le E|X_1|<\infty$,
we have $EL_1^{-1}<\infty$ (e.g., Theorem 1 of \cite{DM4}).
Thus letting $t\to\infty$ in
\begin{equation}\label{normL}
\frac{X_{L_t^{-1}}}{L_t^{-1}}
= \frac{H_t}{L_t^{-1}}
= \left(\frac{H_t}{t}\right)\left(\frac{t}{L_t^{-1}}\right),
 \end{equation}
and using the strong law, we obtain $H_t/t\to EX_1EL_1^{-1}$ as $t\to\infty$.  This implies $EH_1<\infty$ and $EX_1=EH_1/EL_1^{-1}$.
Conversely,  $EH_1<\infty$ implies  $0\le EX_1\le E|X_1|<\infty$
by  Theorem 8 of \cite{DM3}, and $EL_1^{-1}<\infty$
implies  $X$ drifts to $+\infty$ a.s.
by Theorem 1 of \cite{DM4},
so in fact   $0<EX_1\le E|X_1|<\infty$.

(b)(ii) When $X$ is of bounded variation, then $\sigma^2=0$, and by taking limits as $t\downarrow 0$ in  \eqref{normL}, we obtain
\be\label{dlh}
\rmd_{L^{-1}} \rmd_X= \rmd_H.
\ee
When $\rmd_X> 0$, $X_t/t\to \rmd_X>0$ and   $\tau_u/u\to 1/\rmd_X<\infty$ ( by (a)(ii)).  Thus
\ben
\frac{X_{\tau_u}}u = \frac{X_{\tau_u}}{\tau_u}\frac{{\tau_u}}u\to 1\ {\rm\  a.s.\ as }\ u\downarrow 0.
\een
This implies $\rmd_H>0$ by Theorem 4 of \cite{DM3}. Hence by \eqref{dlh} $\rmd_{L^{-1}}>0$ also, and so $c=\rmd_X=\rmd_H/\rmd_{L^{-1}}$.

Conversely assume $\rmd_H>0$, $\rmd_{L^{-1}}>0$ and  $\sigma^2=0$.
We show  $\displaystyle\lim_{u\downarrow  0}{\frac{\tau_u}{u}} =\frac{{\rm d}_{L^{-1}}}{{\rm d}_{H}}$ a.s.
Let $T_u:=\inf\{t>0:H_t>u\}$, $u>0$.  Then $\tau_u=L^{-1}_{T_u}$.  Hence by \eqref{normL}
\ben
\frac{X_{\tau_u}}{u}\frac{u}{\tau_u}
= \left(\frac{H_{T_u}}{{T_u}}\right)\left( \frac{{T_u}}{L_{T_u}^{-1}}\right)\to \frac{\rmd_H}{\rmd_L^{-1}}
\een
Since $\lim_{u\to 0}X_{\tau_u}/u=1$ a.s. when $\rmd_H>0$ by Theorem 4 of \cite{DM3}, the result follows.
\hfill\halmos

\bigskip \noindent {\bf Proof of Theorem \ref{Etau}}.\
We begin by recalling that from Theorem 1 of  \cite{DM4},  $E\tau_u<\infty$ for some, hence all, $u\ge 0$,
iff $X$ drifts to $+\infty$ a.s., iff $EL_1^{-1}<\infty$.

Use identity (8) on p.174 of Bertoin (1996)
to write
\begin{equation}\label{Etid}
E\tau_u=
\lim_{\lambda\dto 0}\frac{\kappa(\lambda,0)}{\lambda}
V_H(u)
= EL_1^{-1} V_H(u),\ \ u>0,
\end{equation}
where
\ben
V_H(u)= \int_0^\infty P(H_t \le u) \rmd t
\een
is the renewal function associated with $H$.

(a)
(i)\ For $u\to\infty$:
by the elementary renewal theorem (Kyprianou  \cite{kypbook}, Cor 5.3 p.114; note
that there is no non-lattice restriction on the support of $\Pi_X$,
and the case $EH_1=\infty$ is covered, e.g., by Gut \cite{Gut04} Theorem 4.1 p.51)
we have
\begin{equation}\label{V1}
\lim_{u\to \infty}\frac{V_H(u)}{u}
=\frac{1}{EH_1}\in [0,\infty),
\end{equation}
so we see that
$\lim_{u\to \infty}E\tau_u/u=1/c$ for some $c\in(0,\infty)$
iff $EH_1<\infty$ and $EL_1^{-1}<\infty$, and then $c= EH_1/EL_1^{-1}$.
Since
$EH_1<\infty$ is equivalent to
$0<EX_1\le E|X_1|<\infty$ when $X_t\to\infty$ by Theorem 8 of \cite{DM3}, we have only left to observe that by Wald's equation for L\'evy processes, \cite{Hall}, $EH_1/EL_1^{-1}=EX_1$.

(ii)\ For $u\dto 0$: assume that $EL_1^{-1}<\infty$ and  $\rmd_H>0$.  Since ${\cal H}_t\ge \rmd_H t$, $t\ge 0$,
it follows easily that
\ben\ba
V_{\cal H}(u)&:= \int_0^\infty P({\cal H}_t \le u) \rmd t\\
&= \int_0^{u/\rmd_H} P({\cal H}_t \le u) \rmd t\\
&\le e^{qu/\rmd_H} \int_0^{u/\rmd_H} e^{-qt}P({\cal H}_t \le u) \rmd t\\
&=e^{qu/\rmd_H}V_H(u),
\ea\een
while trivially
$ V_H(u)\le V_{\cal H}(u)$.
Thus by Theorem III.5 of Bertoin \cite{Bert}, which applies to proper subordinators, we have
\begin{equation}\label{Vlim}
\lim_{u\dto 0}\frac{V_H(u)}{u} =\frac{1}{\rmd_H},
\end{equation}
and so $\lim_{u\dto 0}E\tau_u/u=EL_1^{-1}/\rmd_H$ by \eqref{Etid}.
Conversely,  $\lim_{u\to \infty}E\tau_u/u=1/c$ implies by \eqref{Etid}
that \eqref{Vlim} holds with  $\rmd_H$ replaced by $c EL_1^{-1}>0$.
By Lemma 4, p.52 of \cite{doneystf}, we have
 \begin{equation}\label{Vbd}
\frac{V_H(u)}{u}\asymp \frac{1}{\rmd_H+\int_0^u\pibar_H(y)\rmd y+uq},\
{\rm for\ all}\ u>0.
 \end{equation}
Hence $\rmd_H>0$,
since  $\rmd_H=0$ would imply  $\lim_{u\to \infty}V_H(u)/u=\infty$,
a contradiction.
\hfill\halmos

(b)\
In case $c=0$, we see from \eqref{Etid} and \eqref{V1} that
$\lim_{u\to\infty}E\tau_u/u$ exists and is in $[0,\infty)$ when
$E\tau_u<\infty$  for all $u>0$.
When $\lim_{u\to 0}E\tau_u/u=\infty$,
\eqref{Etid} and \eqref{Vbd} show that $\rmd_H=0$, and conversely.

(c)\
For the case $c=\infty$,
supposing $E\tau_u<\infty$ for each $u>0$,
\eqref{Etid} and \eqref{V1} show that
$\lim_{u\to \infty}E\tau_u/u=0$ iff   $EL_1^{-1}<\infty$ and $EH_1=\infty$,
while for $u\dto 0$, the case $c=\infty$ cannot arise;
$\liminf_{u\to 0}E\tau_u/u>0$
follows from \eqref{Etid} and \eqref{Vbd}.
\hfill\halmos

The following lemma, which is a L\'evy process version of a result of
\cite{lai} for random walk, is needed in the proofs of Theorems
\ref{Grs}
and \ref{jet}.

\begin{lem}[Uniform Integrability of $\tau_u$]\label{Lailem}\
Suppose $X$ is a L\'evy process with
$0<EX_1\le E|X_1|<\infty$
and $\tau_u=\inf\{t>0:X_t>u\}$, $u>0$.
Then $\tau_u/u$ are uniformly integrable as $u\to\infty$, i.e.,
 \begin{equation}\label{ui1}
\lim_{x\to\infty}\limsup_{u\to\infty}
E\left(\frac{\tau_u}{u}{\bf 1}_{\{\frac{\tau_u}{u}>x\}}\right)=0.
 \end{equation}
\end{lem}

\bigskip \noindent {\bf Proof of  Lemma \ref{Lailem}}.\
The  random walk result of \cite{lai} 
can be transferred using a stochastic bound due to Doney
\cite{D2004}. First consider the case when $\Pi\equiv 0$. Then
$X_t=t\gamma+\sigma B_t$, where $\gamma=EX_1>0$, $\sigma\ge 0$, and
$(B_t)_{t\ge 0}$ is a standard Browian motion. In this case $\tau_u$
has an inverse Gaussian distribution and
$E\tau^2_u=u\sigma^2/\gamma^3$, which immediately implies uniform
integrability of $\tau_u/u$.

So, assume $\Pi$ is not identically 0. Then $\pibar(x_0)>0$ for some
$x_0>0$, and by rescaling if necessary we can assume
$c_1:=\pibar(1)>0$. As in  \cite{DM4} , let $\sigma_0=0$ and let
$\sigma_n$, $n=1,2,\ldots$, be the successive times at which $X$
takes a jump of absolute value greater than 1. Then
$e_i:=\sigma_i-\sigma_{i-1}$ are i.i.d exponential rvs with
$E(e_1)=1/c_1$. Define $S_n:=X_{\sigma_n}$,  $n=1,2,\ldots$, and
$\tau^S_u=\min\{n\ge 1:S_n>u\}$, $u>0$. Then $S_n$ is a  random walk
with step distribution $Y_i:= X_{\sigma_i}-X_{\sigma_{i-1}}\eqdr
X_{e_1}$. By Wald's equation,
$EY_1=EX_1/c_1>0$. Thus by \cite{lai},  $\tau^S_u/u$ are uniformly
integrable.

Now we can use similar calculations as on p.287 of  \cite{DM4} to
bound the expression on the left of \eqref{ui1} in terms of a
similar expression involving  $\tau^S_u$.  For any $Z\ge 0$ and any
$a>0$ \be\label{Zt} E(Z;Z>a)=\int_{z>a}P(Z>z) \rmd z + aP(Z>a). \ee
Taking $u>0$, $x>0$, $c>1/c_1$, we obtain \be\label{ui2}
E\left(\frac{\tau_u}{uc}{\bf 1}_{\{\frac{\tau_u}{uc}>x\}}\right) =
u^{-1}\int_{y>xu}P\left(\tau_u>yc\right)\rmd y
+xP\left(\tau_u>xuc\right). \ee By Theorem \ref{taurs} we have
$\tau_u/u\topr 1/EX_1$ as $u\to\infty$. So the second term on the
righthand side of \eqref{ui2} tends to 0 as $u\to\infty$ once
$xc>1/EX_1$. {}As on  p.287 of  \cite{DM4}  we have
\[
P\left(\sigma_j\le\tau_u< \sigma_{j+1}\right) =P\left(\widetilde
m_0\le u, \tau^S_{u-\widetilde m_0}=j\right),\ j\ge 1,
\]
where $\widetilde m_0$ is a finite rv independent of
$(S_n)_{n=1,2,\ldots}$. The first term on the righthand side of
\eqref{ui2} is bounded by
\begin{equation}\label{ui3}
u^{-1}\sum_{n\ge\lfloor xu\rfloor} P(\tau_u>nc),
\end{equation}
and now we argue as follows:
 \begin{eqnarray}\label{ui4}
 u^{-1}\sum_{n\ge\lfloor xu\rfloor} P(\tau_u>nc)
&\le&
 u^{-1}\sum_{n\ge\lfloor xu\rfloor} P(\tau_u\ge\sigma_n)
+ u^{-1}\sum_{n\ge\lfloor xu\rfloor} P(\sigma_n>nc)
 \nonumber\\ &=&
u^{-1}\sum_{j\ge\lfloor xu\rfloor}(j-\lfloor xu\rfloor +1)
P\left(\sigma_j\le\tau_u< \sigma_{j+1}\right)
  \nonumber\\
&&\qquad\qquad\qquad
 + u^{-1}\sum_{n\ge1} P\left(\sum_{i=1}^n e_i>nc\right).
 \end{eqnarray}
Since $e_i$ are i.i.d. with a finite exponential moment and
$Ee_i=1/c_1<c$, the sum in the second term on the righthand side of
\eqref{ui4} is convergent, and hence this  term is $o(1)$ as
$u\to\infty$.
 The first term on the righthand side of \eqref{ui4} is
\begin{eqnarray*}
&&u^{-1}\sum_{j\ge\lfloor xu\rfloor}(j-\lfloor xu\rfloor +1)
P\left(\widetilde m_0\le u, \tau^S_{u-\widetilde m_0}=j\right)
\nonumber\\ &=& u^{-1}\sum_{j\ge\lfloor xu\rfloor}(j-\lfloor
xu\rfloor+1) \int_0^uP\left(\tau^S_{u-y}=j\right)P(\widetilde
m_0\in\rmd y) \nonumber\\ &=& u^{-1}\sum_{n\ge\lfloor xu\rfloor}
\int_0^uP\left(\tau^S_{u-y}\ge n\right) P(\widetilde m_0\in\rmd y)
 \nonumber\\ &\le&
 u^{-1}\sum_{n\ge \lfloor xu\rfloor}P\left(\tau^S_u\ge n\right)
  \nonumber\\ &\le&
u^{-1}\int_{y>\lfloor xu\rfloor-1}P\left(\tau^S_u>y\right)\rmd y
 \nonumber\\ &\le&
E\left(\frac{\tau^S_u}{u}{\bf 1}_{\{\frac{\tau^S_u}{u}> \frac
x2\}}\right),
\end{eqnarray*}
if $xu\ge 4$, where the last inequality comes from \eqref{Zt}. Since
$\tau^S_u/u$ are uniformly integrable by Lai's result we get
\eqref{ui1}.
\hfill\halmos

\begin{rem}  Lemma \ref{Lailem} could  be used to give an alternative proof of  Part (a)(i) of Theorem
\ref{Etau} from Theorem \ref{taursas}.  This approach however would not work for Part (a)(ii) of Theorem
\ref{Etau}, since $\tau_u/u$ are not uniformly integrable as $u\downarrow 0$.  This is because the almost sure limit in Theorem \ref{taursas} does not agree with the limit in mean in Theorem \ref{Etau} as $u\downarrow 0$.
\end{rem}

\bigskip \noindent {\bf Proof of Theorem \ref{Grs}}.\
(a) We first observe that \eqref{rt3} is equivalent to
\begin{equation}\label{rt3a}
\lim_{x\to 1/L}\frac{\kappa(0,x)-\kappa(0,0)}{\kappa(\xi x,x)}=\frac{c}{c+\xi},
\  {\rm for\ each}\ \xi>0.
\end{equation}
When $L=\infty$ this because $\kappa(0,0)=q=0$ by our assumption that $\limsup_{t\to\infty}X_t=\infty$.  When  $L=0$ we have, by \eqref{kapexp}, that $\kappa(\xi x,x)\to\infty$ as $x\to 1/L$ for each $\xi>0$ unless $\rmd_H=\rmd_{L^{-1}}=0$ and $\Pi_{L^{-1},H}$ is a finite measure. But this is impossible since $0$ is regular for $(0,\infty)$, so \eqref{rt3a} is equivalent when $L=0$ also.

From \eqref{e3} we have, for $u>0$, $\nu>0$,
\begin{equation*}
\mu\int_{u\ge 0}  e^{-\mu u}
E\left(e^{-\nu G_{\tau_u-}};\tau_{u}<\infty\right)\rmd u
=\frac{\kappa(0,\mu)-\kappa(0,0)}{\kappa(\nu,\mu)}.
\end{equation*}
Take $y>0$ and replace $\mu$ by $\mu/y$, $u$ by $uy$, and $\nu$
by $\nu/y$ in this to get
\begin{equation}\label{rsm}
\mu\int_{u\ge 0}  e^{-\mu u}
E\left(e^{-\nu G_{\tau_{uy-}}/y};\tau_{uy}<\infty\right)\rmd u
= \frac{\kappa(0,\mu/y)-\kappa(0,0)}{\kappa(\nu/y,\mu/y)}.
 \end{equation}
\eqref{rt4} implies $ G_{\tau_{uy-}}/y\topr u/c$, as $y\to L$,
for each $u>0$. So the lefthand side of \eqref{rsm} tends to
$c/(c+\nu/\mu)$ as $y\to L$, and then \eqref{rt3} follows
from \eqref{rt3a} and  the righthand side of \eqref{rsm}.
The proof that  \eqref{rt3} implies
\eqref{rt4}  is analogous to that  in Theorem \ref{taurs}.

(b) By the strong law when $L=\infty$ and by Theorem 39 of \cite{doneystf} when $L=0$, it suffices to prove that for $c\in [0,\infty)$
\be\label{X=c}
\lim_{t\to L}{\frac{X_t}{t}} = c\ \rm{a.s.},
\ee
if and only if
\be\label{G=c}
\lim_{u\to L}{\frac{G_{\tau_u-}}{u}} = \frac{1}{c}\ \rm{a.s.}
\ee

A simple pathwise argument shows that  \eqref{X=c} implies  \eqref{G=c} when $c>0$, but the case $c=0$ is a little trickier.  Since the following argument works whenever $c\in[0,\infty)$, we prove it under that assumption.  So assume that \eqref{X=c} holds with $c\in[0,\infty)$.  Then $\Xbar_t/t\to c$ as $t\to L$.
If \eqref{G=c} fails, then
\ben
\liminf_{u\to L}{\frac{G_{\tau_u-}}{u}} < \frac{1}{c},
\een
since under \eqref{X=c}, $\displaystyle \frac{G_{\tau_u-}}{u}\le \frac{{\tau_u}}{u}\to \frac 1c$ as $u\to L$, by Theorem \ref{taursas}.
Now consider the random level $Z_u=\Xbar_{\tau_u-}$, and observe that
\ben
\tau_{Z_u}=\tau_u, \ \   \Xbar_{G_{\tau_u-}}=\Xbar_{{\tau_u-}}=Z_u.
\een
Writing
\ben
\frac{\Xbar_{G_{\tau_u-}}}{G_{\tau_u-}}=\frac{Z_u}{u}\frac{u}{G_{\tau_u-}},
\een
it then follows that
\ben
\liminf_{u\to L}\frac{Z_u}{u}<1.
\een
Thus
\ben
\limsup_{u\to L}\frac{X_{\tau_{Z_u}}}{Z_u}=\limsup_{u\to L}\frac{X_{\tau_u}}{Z_u}\ge \limsup_{u\to L}\frac{u}{Z_u}>1.
\een
In particular it is not the case that $\displaystyle \lim_{v\to L}\frac{X_{\tau_{v}}}{v}= 1$ a.s.  This implies $E|X|=\infty$ when $L=\infty$, by Theorem 8 of \cite{DM3}, and $X$ is not of bounded variation when $L=0$, by Theorem 4 of \cite{DM3}. In either case \eqref{X=c} fails to hold which is a contradiction.

Conversely assume \eqref{G=c} holds for some $c\in [0,\infty)$.  Then  for any $a>c$,
$\displaystyle{G_{\tau_u-}}> {u}{a}^{-1}$ eventually.  Hence $\displaystyle\Xbar_{{u}{a}^{-1}}\le u$ eventually.  This implies
\be\label{lec}
\limsup_{t\to L} \frac{X_t}{t}\le c.
\ee
If $L=\infty$, then arguing as in the proof of Theorem \ref{taursas} (a)(i),  \eqref{lec}
implies that $0\le EX_1\le E|X_1|<\infty$ and so
\eqref{X=c} holds with $c=EX_1\ge 0$. Since \eqref{X=c} implies  \eqref{G=c},  the constant $c$ for which \eqref{G=c} was assumed to hold must also have been $c=EX_1$, completing the proof of (i).  If $L=0$ then \eqref{lec} forces $X$ to have bounded variation with
$\rmd_X\ge 0$ since $0$ is regular for $(0,\infty)$.  In that case  \eqref{X=c} holds with $c=\rmd_X$, and so again since \eqref{X=c} implies  \eqref{G=c}, the constant $c$ for which \eqref{G=c} was assumed to hold must also have been $c=\rmd_X$, completing the proof of (ii).

(c) Finally, suppose $0<EX_1\le E|X_1|<\infty$. Then \eqref{G=c} holds with $c=EX_1$ and, further, $\Gtau/u$ are uniformly integrable as $u\to\infty$
by Lemma \ref{Lailem}.
Thus we get $\lim_{u\to\infty}E(\Gtau/u)=1/EX_1$.
\hfill\halmos\bigskip

\bigskip \noindent {\bf Proof of Theorem \ref{jet}}.\
Since $0<EX_1\le E|X_1|<\infty$, we have by Theorems \ref{taursas} and \ref{Grs}
\ben
\lim_{u\to \infty}{\frac{G_{\tau_u-}}{u}} = \lim_{u\to \infty}{\frac{{\tau_u}}{u}}=\frac{1}{EX_1}\ {\rm a.s.}
\een
The assumptions on $X$ imply that $H$ does not have a lattice jump
distribution, and hence it follows from
 \cite{BVS}
that
\[
X_{\tau_u}-u \todr Y, \  {\rm  as}\  u\to\infty,
\]
where $Y$ is the rv defined in the statement of Theorem \ref{jet}.
Since $\tau_u/u$, and consequently $\Gtau/u$ also, are uniformly integrable as $u\to\infty$, by Lemma \ref{Lailem}, the result follows.
\hfill\halmos \bigskip

\setcounter{equation}{0} \section{Proofs for Section
\ref{s3}}\label{s6}

We assume throughout this section the setup of Section \ref{s3}.
Let ${\cal F}_{\tau_u}$ be the $\sigma$-algebra generated by $X$ up to time $\tau_u$. By
Corollary 3.11 of \cite{kypbook}, for any $Z_u$ which is nonnegative and
measurable with respect to ${\cal F}_{\tau_u}$, we have
 \begin{equation}\label{*form}
E\left(Z_u;\tau_u<\infty\right)
= E^*\left(Z_ue^{-\nu_0 X_{\tau_u}}\right).
 \end{equation}
This immediately yields
the following lemma,
 which can be found  in Theorem IV.7.1 of \cite{aa} for
 compound Poisson processes with negative drift.
 Our proof is analogous to that  in \cite{aa}.

\begin{lem}\label{As2.4}
Suppose $\mu^*<\infty$ and $Y_u$ are
${\cal F}_{\tau_u}$-measurable rvs such that
$Y_u\toprstar 0$ as $u\to\infty$. Then $Y_u\topru 0$.
\end{lem}

\bigskip \noindent {\bf Proof of Lemma \ref{As2.4}}:\
For $\veps>0$, by \eqref{*form},
 \begin{eqnarray*}
P^{(u)}(|Y_u|>\veps)
= \frac{P(|Y_u|>\veps, \tau_u<\infty)}{P( \tau_u<\infty)}
=\frac{ E^*\left(e^{-\nu_0( X_{\tau_u}-u)};|Y_u|>\veps\right)}
{e^{\nu_0 u}P( \tau_u<\infty)}.
\end{eqnarray*}
Since  $\mu^*<\infty$ we have by \eqref{BDmain} that
the denominator here is bounded away from 0, hence the result.
\hfill \halmos

\bigskip \noindent {\bf Proof of Theorem \ref{thm1}}\
Since $X$ is a L\'evy process under $P^*$ with $E^*X_1=\mu^*\in(0,\infty)$,
it follows easily from the strong law that
\be\label{*cvg}
\frac{X_{\tau_u}}{u}\to 1,\ \frac{G_{\tau_u-}}{u}\to \frac{1}{\mu^*}
\quad {\rm and}\quad
\frac{\tau_u}{u}\to \frac{1}{\mu^*},
\ P^*-a.s., \  {\rm as}\ u\to \infty.
 \ee
\eqref{Astau} is then immediate from  Corollary \ref{As2.4}.

For \eqref{Eu0}, use Theorem \ref{jet} , \eqref{BDmain} and  \eqref{*form}
to deduce that, as $u\to\infty$,
 \begin{eqnarray}\label{Eu}
\frac{E^{(u)}\tau_u}{u}
=\frac{E\left(\tau_u;\tau_u<\infty\right)}{uP(\tau_u<\infty)}
=\frac{E^*\left(\tau_u e^{-\nu_0(X_{\tau_u}-u)}\right)}
{ue^{\nu_0 u}P(\tau_u<\infty)}
\to \frac{1}{C\mu^*}E^*e^{-\nu_0 Y}
=\frac{1}{\mu^*}.
\end{eqnarray}
The limit involving $G_{\tau_u-}$ is similar.
For the final  limit in \eqref{Eu0}, first observe that
\[
 \frac{E^{(u)}X_{\tau_u}}{u}
= \frac{E\left(X_{\tau_u};\tau_u<\infty\right)}{uP(\tau_u<\infty)}
\sim \frac{E^*\left(X_{\tau_u} e^{-\nu_0(X_{\tau_u}-u)}\right)}{Cu}.
\]
Now let $O_u:=X_{\tau_u}-u$, $u>0$.
Then $\left(X_{\tau_u}/u\right)e^{-\nu_0 O_u}$ is
uniformly integrable, because for $x>1$ and
$c_0:=\sup_{y\ge 0}(ye^{-\nu_0 y})=(e\nu_0)^{-1}$ we have
 \begin{eqnarray*}
E^*\left(\frac{X_{\tau_u}}{u}
e^{-\nu_0 O_u}{\bf 1}_{\{\frac{X_{\tau_u}}{u}>x\}}\right)
&=&
 E^*\left(\frac{O_u}{u} e^{-\nu_0 O_u}
{\bf 1}_{\{O_u>(x-1)u\}}\right)
 +E^*\left(e^{-\nu_0 O_u} {\bf 1}_{\{O_u>(x-1)u\}}\right)
\nonumber\\
&\le&
\frac{c_0}{u}+e^{-\nu_0(x-1)u}.
\end{eqnarray*}
Letting $u\to\infty$ then $x\to\infty$ shows the  uniform
integrability.
Since
$X_{\tau_u}-u\todrstar Y$
and  $X_{\tau_u}/u\toprstar 1$ by \eqref{*cvg}, we have
\[
 \frac{E^{(u)}X_{\tau_u}}{u}\sim C^{-1}E^*\left(\frac{X_{\tau_u}}{u} e^{-\nu_0 O_u}\right)
\to C^{-1}E^*e^{-\nu_0 Y} =1,
\]
completing the proof. \hfill \halmos

\bigskip \noindent {\bf APPENDIX}\

\bigskip \noindent {\bf Proof of Lemma \ref{cex}}.\
We first construct a L\'evy process satisfying \eqref{ce1}.  For the characteristics of $X$ we take $\gamma=-2$, $\sigma=0$ and the L\'evy measure given by
\ben
\pibar^+_X(x)=\frac 1{x|\ln x|},\ \ \pibar^-_X(x)=\pibar^+_X(x)+\frac {\ln 2}{x(\ln x)^2},\ \ 0<x<1/2
\een
\ben
\pibar^+_X(x)=\pibar^-_X(x)=0,\ \ x\ge 1/2.
\een
Then $X$ is not of bounded variation since
\ben
\int (|x|\wedge 1)\Pi(\rmd x)=\infty,
\een
and consequently $0$ is regular for $(0,\infty)$.
Further one can easily check that \eqref{rt46} holds with $c=-1$, and so
\ben
\frac{X_t}{t}\topr -1\ {\rm as}\ t\dto 0
\een
by Theorem 2.1 of \cite{DM2}.  Also, the argument given in Theorem 2.1 that \eqref{rt5} implies \eqref{rt41} when $L=0$, shows that under \eqref{rt46},
\ben
P(\sup_{0\le s\le t}|X_s+s|>\varepsilon t)\to 0\ {\rm as}\ t\dto 0.
\een
From this we conclude that
\ben
\frac{\Xbar_t}{t}\topr 0\ {\rm as}\ t\dto 0,
\een
completing the example.

\medskip
We now construct a L\'evy process satisfying \eqref{ce}.  This is based on  Example 3.5 in \cite{KM99}, which constructs a random walk $S_n=\sum_{i=1}^nY_i$
with i.i.d. summands $Y_i$ which satisfying
\begin{equation}\label{Slim}
\frac{S_n }{ n} {\buildrel P \over \longrightarrow} - \infty \quad
{\rm and}\quad
{\Sbar_n \over n} = {\max_{1\le j\le n}S_j \over n}
{\buildrel P \over \longrightarrow} \infty.
 \end{equation}
This is done by finding a random walk which
is negatively relatively stable (NRS) as $n\to \infty$, i.e.,
with
\begin{equation}\label{SD}
{S_n \over D(n)} {\buildrel P \over \longrightarrow} - 1, \  (n
\rightarrow \infty),
\end{equation}
for a norming sequence $D(n)>0$ with  $D(n)/n\to\infty$  satisfying
\begin{equation}\label{YD}
\lim_{n\to\infty}
\sum^n_{j = \ell(n)} P(S_1> x D (j))= \infty,
 \end{equation}
where $\ell(n)$ is inverse to $D(n)$. The sequences $D(n)$ and
$\ell(n)$ are
strictly increasing to $\infty$ as $n\to \infty$ and satisfy
$D(n)=-nA(D(n))$
(where $A(\cdot)$ is defined in \eqref{Adef})
and $\ell(n)=-n/A(n)$. The function $-A(x)$ is positive for $x$
large enough, slowly varying as $x\to \infty$, and tends to $\infty$ as $x\to \infty$.

(i)\
Consider first the case $L=\infty$.
Let $(N_t)_{t\ge 0}$ be a Poisson process of rate $1$ independent of
the $Y_i$ and
set $X_t:= S_{N_t}$,  $t\ge 0$, where $S_n$ is as in \eqref{Slim}.
Then the compound Poisson process $X_t$ satisfies \eqref{ce} with $
L=\infty$. This is fairly straightforward to check and we omit the
details.

(ii)\
Now consider the case $L=0$.
For this we have to modify  Example 3.5 in \cite{KM99} to work as
$t\dto 0$.
Details are as follows.

We construct a L\'evy process $X_t$ which is NRS
as $t\dto 0$, i.e., is such that
\begin{equation}\label{NRS}
\frac{X_t}{b(t)} {\buildrel P \over \longrightarrow} -1
\ {\rm as}\ t\dto 0,
 \end{equation}
for a nonstochastic function $b(t) > 0$, with $b(t)\dto 0$ and
$b(t)/t\to\infty$ as $t\dto 0$.
To do this, it will be useful to summarize here some properties
concerning
(negative) relative stability at 0 of $X_t$;
for reference, see \cite{DM2} (and replace $X$ by $-X$).
We assume that $\pibar_X(x) > 0$  for all $x>0$.
We then have that $X_t \in NRS$
if and only if $\sigma=0$ and the function $A(x)$
defined in \eqref{Adef} is strictly negative for
all $x$ small enough, $x \le  x_0$, say, and satisfies
\begin{equation}\label{ap}
\lim_{x \dto 0} {A(x) \over x \pibar_X(x)} =-\infty.
\end{equation}

When $A(x) < 0$ for $x \le x_0$, there is an $x_1\le x_0$
so that the function
 \begin{equation}\label{1.15}
D(x): = \sup \left\{y \ge x_0: {-A (y) \over y} \ge {1\over
x}\right\},\ 0<x\le x_1,
\end{equation}
is strictly positive and finite and satisfies
$$
D(x) =- xA(D(x))
$$
for all $x \le  x_1$.
It is easily seen to be strictly
increasing on $x \le x_1$ (by the continuity of $y \mapsto
-A(y)/y$) with $D(0)=0$. Thus for small enough $y$
we can define the inverse function
$$
\ell (y) = \sup \{x: D(x) \le y\} = \inf\{x:D(x) > y\}.
$$
When \eqref{ap} holds, $-A(x)$ is slowly varying as $x\dto 0$,
and as a consequence $D(x)$ and $\ell(x)$
are both regularly varying with  index 1 as $x \dto 0$
(see Bingham, Goldie and Teugels, \cite{BGT}, Theorem 1.5.12).
Also, when \eqref{ap} holds, it is easy to check that
$-A(y)/y$ is strictly decreasing for $y$ small enough,
so $\ell(y)$ is continuous and strictly increasing, with $\ell(0)=0$,
 and
$$
\ell(y) = \frac y{-A(y)},
$$
for $y$ small enough.
Finally, we can take $b(t) = D(t)$ in \eqref{NRS} to get
negative relative stability of $X$ in the form:
 \begin{equation}\label{5.29}
\frac{X_t}{D(t)} {\buildrel P \over \longrightarrow} -1
\ {\rm as}\ t\dto 0.
 \end{equation}

To construct the process required in the lemma,
we will specify a L\'evy measure
$\Pi_X$ for $X$ such that
\begin{equation}\label{5.29a}
\frac{D(t)} t \rightarrow \infty, \ {\rm as}\ t\dto 0,
\end{equation}
consequently $\ell(t)/t\to 0$, and
\begin{equation}\label{5.30}
\int_{\ell(t)}^t \pibar^+_X(xD(s))\rmd s \to \infty,
\ {\rm as}\ t\dto 0,
\end{equation}
for all $x  > 0$.
\eqref{5.29} then implies $X_t/t {\buildrel P \over \longrightarrow} -
\infty$
as $t\dto 0$. We claim that in addition \eqref{5.29}--\eqref{5.30}
imply
\begin{equation}\label{5.31}
{\Xbar_t \over t} {\buildrel P \over \longrightarrow} \ \infty\ {\rm
as}\ t\dto 0.
\end{equation}

To prove \eqref{5.31} from  \eqref{5.29}--\eqref{5.30},
fix $t>0$ and $x > 1$ and write
\begin{eqnarray}\label{5.32}
&&
P\left(\Xbar_t  > x t\right)
=
\lim_{n\to\infty}
P\left(\max_{1\le j\le nt} X(j/n)  > x D(\ell(t)) \right)
\nonumber\\
&\ge&
\liminf_{n\to\infty}
\sum_{n\ell(t)\le j\le nt}
P\left(X((j-1)/n)>-x D(j/n),\
\max_{j+1\le k\le nt}
\frac{\Delta(k/n)}{D(k/n)} \le 2x < \frac{\Delta(j/n)}{D(j/n)}
\right),
\nonumber\\&&
\end{eqnarray}
where
\[
\Delta(k/n):= X(k/n)-X((k-1)/n), \ k=1,2, \ldots,
\]
are i.i.d. with distribution the same as that of $X(1/n)$.
Given $\veps\in(0,1)$, by \eqref{5.29} there is a $t_0>0$ such  that
$P\big(X((j-1)/n)>-x D(j/n)\big)>1-\veps$ when $j/n\le t\le t_0$.
Thus, keeping $t\le t_0$,
\begin{eqnarray*}
P\left(\Xbar_t  > x t\right)
&\ge&
(1-\veps) \liminf_{n\to\infty}
\sum_{n\ell(t)\le j \le nt}
P\left( \max_{j+1\leq k\leq nt} \frac{\Delta(k/n)}{D(k/n)} \le 2x
< \frac{\Delta(j/n)}{D(j/n)} \right).
\end{eqnarray*}
Here the sum equals
\begin{eqnarray*}
&&P\left(\Delta(j/n)>2x D(j/n),\ {\rm for\ some}\
j\in[n\ell(t),nt]\right)
=
1-\prod_{n\ell(t)\le j\le nt} P\big(X(1/n)\le 2x D(j/n)\big)
\nonumber\\
&& \qquad \qquad \qquad
\ge
1-\exp\left(-\sum_{n\ell(t)\le j\le nt}P\big(X(1/n)>2x D(j/n)\big)\right).
\end{eqnarray*}
Noting that
\begin{eqnarray*}
n\int_{\ell(t)}^{t+1/n} P\big(X(1/n)> 2xD(s)\big)\rmd s
&\le&
 n\sum_{n\ell(t)\le j\le nt}\int_{j/n}^{(j+1)/n} P\big(X(1/n)> 2xD(s)\big)\rmd
s
\nonumber\\ &\le&
\sum_{n\ell(t)\le j\le nt}P\big(X(1/n)> 2xD(j/n)\big)\rmd s,
\end{eqnarray*}
and employing the fact that
$\lim_{n\to\infty}nP(X(1/n)>a)=\pibar_X^+(a)$ for all $a>0$
(cf. \cite{Bert} p.39) we get
\begin{eqnarray*}
P\left(\Xbar_t  > x t\right)
&\ge&
(1-\veps) \liminf_{n\to\infty}
\left(1-\exp\left(-n \int_{\ell(t)}^{t+1/n}P\big(X(1/n)> 2xD(s)\big)\rmd
s\right)\right)
 \nonumber\\
&=&
(1-\veps)
\left(1-\exp\left(-\int_{\ell(t)}^t \pibar^+_X(2xD(s)\rmd
s\right)\right).
\end{eqnarray*}
The last expression
 tends to 1 as $t\dto 0$ then $\veps\dto 0$, provided \eqref{5.30}
holds.
Thus \eqref{5.31} will follow from \eqref{5.29}--\eqref{5.30}, as
claimed.

It remains to give an example where \eqref{5.29}--\eqref{5.30} hold.
Define
$$
L(x)=e^{(-\log x)^\beta}, \ 0<x < e^{-1},
$$
and keep $ {1\over 2} < \beta < 1$.  Choose a  L\'evy measure $\Pi_X$
which satisfies
$$
\pibar^+_X(x) = -2L^\prime (x) = {2\beta (-\log x)^{\beta - 1} L
(x)\over x}
$$
and
$$
\pibar^-_X(x)= -L'(x) =\pibar^+_X(x)/2,
$$
for $x$ small enough, $x \le  x_0$, say.
(Note that $\pibar^+_X(x)$ and $\pibar^-_X(x)$ are infinite at 0
and decrease to 0 as $x \to \infty$.)   A straightforward calculation
using \eqref{Adef}, gives,  for $x > 0$,
\begin{eqnarray*}
A(x) &=&
\gamma +\pibar_X^+(1)- \pibar_X^-(1)+1-L (x).
\end{eqnarray*}
Thus $A(x) \rightarrow  - \infty$ as $x\dto 0$, and
\begin{equation}\label{5.33}
{-A (x) \over x \pibar_X(x)}
\sim
\frac{L(x)}{-3L'(x)} \sim {(-\log x)^{1- \beta} \over 3\beta}
 \rightarrow \infty \  \hbox{ as} \ x\dto 0.
\end{equation}
Thus by \eqref{ap}, $X$ is negatively relatively stable at 0 and
\eqref{5.29} holds with $D(t)$ as in \eqref{1.15}.
We have $D(t) =- t A (D(t))$ and $\ell (t) =-t/A(t)$
for small enough $t>0$.
Since $-A(x) \to \infty$ as $x \dto 0$, \eqref{5.29a} also holds.
Now $L(x)$ is slowly varying as $x\dto 0$,
so $\pibar^+_X(x)$ is regularly varying with index $-1$
as $x \dto 0$. It therefore suffices
to check \eqref{5.30} with $x = 1$.
Note also that $\log (D (t)) \sim \log t$ as $t\dto 0$, because
$A(\cdot)$ is slowly varying at 0.  Hence
for some constant $c_2 > 0$, as $t\dto 0$,
\begin{eqnarray*}
\int_{\ell(t)}^t \pibar^+_X(D(s))\rmd s
&=&
\int_{\ell(t)}^t
{D(s) \pibar^+_X(D (s)) \over -s A (D(s))}\rmd s
 \nonumber\\ & \sim&
3\beta\int_{\ell(t)}^t
{(-\log D(s))^{\beta - 1} \over s}\rmd s
 \nonumber\\
& \ge &
3\beta(-\log D(t))^{\beta - 1} \left(\log t - \log \ell(t)\right)
 \nonumber\\
\nonumber\\ &\ge&
c_2 (-\log t)^{2\beta - 1}\to \infty
\end{eqnarray*}
(because $\beta > {1\over 2}$).  Thus \eqref{5.30} holds too.
\hfill\halmos

\end{document}